\title{The isoperimetric constant of the random graph process}
\author{{Itai Benjamini \thanks{Weizmann Institute, Rehovot, 76100,
Israel. Email: itai.benjamini@weizmann.ac.il}} \quad {Simi Haber
\thanks{ Department of Mathematics, Raymond and Beverly
Sackler Faculty of Exact Sciences, Tel Aviv University, Tel Aviv,
69978, Israel. Email: habbersi@tau.ac.il.}} \quad {Michael
Krivelevich\thanks{Department of Mathematics, Raymond and Beverly
Sackler Faculty of Exact Sciences, Tel Aviv University, Tel Aviv,
69978, Israel. Email: krivelev@post.tau.ac.il. Research supported in
part by a USA-Israeli BSF grant and a grant from the Israeli Science
Foundation.}} \quad {Eyal Lubetzky
\thanks{ Department of Computer Science, Raymond and Beverly
Sackler Faculty of Exact Sciences, Tel Aviv University, Tel Aviv,
69978, Israel. Email: lubetzky@tau.ac.il.}}}

\documentclass [11pt]{article}
\usepackage[]{epsf,epsfig,amsmath,amssymb,amsfonts,latexsym,amsthm}

\newtheorem{theorem}{Theorem}[section]
\newtheorem{lemma}[theorem]{Lemma}
\newtheorem{claim}[theorem]{Claim}
\newtheorem*{definition}{Definition}

\newtheorem{proposition}[theorem]{Proposition}

\renewcommand{\epsilon}{\varepsilon}

\newcommand{\Small}{\ensuremath{\text{\sc Small}}}
\newtheoremstyle{upright}%
        {8pt plus2pt minus4pt}%
        {8pt plus2pt minus4pt}%
        {\upshape}%
        {}%
        {\bfseries\scshape}%
        {:}%
        {1em}%
        {}%
\theoremstyle{upright}
\newtheorem{remark}[theorem]{Remark}

\begin{document}
\maketitle

\begin{abstract}
The isoperimetric constant of a graph $G$ on $n$ vertices, $i(G)$,
is the minimum of $\frac{|\partial S|}{|S|}$, taken over all
nonempty subsets $S\subset V(G)$ of size at most $n/2$, where
$\partial S$ denotes the set of edges with precisely one end in $S$.
A random graph process on $n$ vertices, $\widetilde{G}(t)$, is a
sequence of $\binom{n}{2}$ graphs, where $\widetilde{G}(0)$ is the
edgeless graph on $n$ vertices, and $\widetilde{G}(t)$ is the result
of adding an edge to $\widetilde{G}(t-1)$, uniformly distributed
over all the missing edges. We show that in almost every graph
process $i(\widetilde{G}(t))$ equals the minimal degree of
$\widetilde{G}(t)$ as long as the minimal degree is $o(\log n)$.
Furthermore, we show that this result is essentially best possible,
by demonstrating that along the period in which the minimum degree
is typically $\Theta(\log n)$, the ratio between the isoperimetric
constant and the minimum degree falls from $1$ to $\frac{1}{2}$, its
final value.
\end{abstract}

\section{Introduction}

Let $G=(V,E)$ be a graph. For each subset of its vertices, $S
\subseteq V$, we define its edge boundary, $\partial S$, as the set
of all edges with exactly one endpoint in $S$: $$\partial S = \{
(u,v) \in E : u \in S, v \notin S\} ~.$$ The isoperimetric constant,
or isoperimetric number, of $G=(V,E)$, $i(G)$, is defined to be:
$$i(G) = \min_{\emptyset \neq S \subset V} \frac{|\partial S|}{\min\{|S|,|V \setminus S|\}}
 = \mathop{\min_{\emptyset \neq S \subset V}}_{|S| \leq \frac{1}{2}|V|} \frac{|\partial S|}{|S|} ~.$$

It is well known that this parameter, which measures edge expansion
properties of a graph $G$, is strongly related to the spectral
properties of $G$, and namely:
\begin{equation}\label{eigenvalue-bounds}
\frac{\lambda}{2} \leq i(G) \leq \sqrt{\lambda(2\Delta(G) -
\lambda)} ~, \end{equation} where $\Delta(G)$ denotes the maximal
degree of $G$, and $\lambda$ denotes the second smallest eigenvalue
of the Laplacian matrix of $G$ (for proofs of these facts, see
\cite{AlonLowerEigenvalueBound} and \cite{MoharIso}). The upper
bound in \eqref{eigenvalue-bounds} can be viewed as a discrete
version of the Cheeger inequality bounding the first eigenvalue of a
Riemannian manifold, and indeed, there is a natural relation between
the study of isoperimetric inequalities of graphs and the study of
Cheeger constants in spectral geometry. For instance, see
\cite{Buser}, where the author relates between isoperimetric
constants and spectral properties of graphs and those of certain
Riemann surfaces. The eigenvalue bounds in \eqref{eigenvalue-bounds}
also relate $i(G)$ (as well as a variation of it, the conductance of
$G$) to the mixing time of a random walk in $G$, defined to be the
minimal time it takes a random walk on $G$ to approach the
stationary distribution within a variation distance of $1/2$.

A closely coupled variant of the isoperimetric constant is the
Cheeger constant of a graph, where the edge boundary of $S$ is
divided by its volume (defined to be the sum of its degrees) instead
of by its size. For further information on this parameter, its
relation to the isoperimetric constant, and its corresponding
eigenvalue bounds (analogous to \eqref{eigenvalue-bounds}), see
\cite{ChungCheeger}, as well as \cite{SpectralGraphTheory}, Chapter
2.

There has been much study on the isoperimetric constants of various
graphs, such as grid graphs, torus graphs, the $n$-cube, and more
generally, cartesian products of graphs. See, for instance,
\cite{BollobasLeaderGrid,BollobasLeaderTorus,ChungTetali,HoudreTetali,
MoharIso}. In \cite{BollobasLowerRegulerBound}, Bollob\'{a}s studied
the isoperimetric constant of random $d$-regular graphs, and used
probabilistic arguments to prove that infinitely many $d$-regular
graphs $G$ satisfy $i(G) \geq \frac{d}{2} - O(\sqrt{d})$. Alon
proved in \cite{AlonUpperRegulerBound} that this inequality is in
fact tight, by providing an upper bound of $i(G) \leq \frac{d}{2} +
O(\sqrt{d})$ for any $d$-regular graph $G$ on a sufficiently large
number of vertices.

In this paper, we study the isoperimetric constant of general random
graphs $\mathcal{G}(n,p)$, $\mathcal{G}(n,M)$, and the random graph
process, and show that in these graphs, the ratio between the
isoperimetric constant and the minimal degree exhibits an
interesting behavior.

We briefly recall several elementary details on these models (for
further information, c.f., e.g., \cite{RandomGraphs}, Chapter 2).
The random graph $\mathcal{G}(n,p)$ is a graph on $n$ vertices,
where each pair of distinct vertices is adjacent with probability
$p$, and independently of all other pairs of vertices. The
distribution of $\mathcal{G}(n,p)$ is closely coupled with that of
$\mathcal{G}(n,M)$, a uniform distribution on all graphs on $n$
vertices with precisely $M$ edges, if we choose $p=M/\binom{n}{2}$.
The random graph process on $n$ vertices, $\widetilde{G}(t)$, is a
sequence of $\binom{n}{2}$ graphs, where $\widetilde{G}(0)$ is the
edgeless graph on $n$ vertices, and $\widetilde{G}(t)$ is the result
of adding an edge to $\widetilde{G}(t-1)$, uniformly distributed
over all the missing edges. Notice that at a given time $0 \leq t
\leq \binom{n}{2}$, $\widetilde{G}(t)$ is distributed as
$\mathcal{G}(n,M)$ with $M=t$.

For a given graph process $\widetilde{G}$ on $n$ vertices, we define
the hitting time of a monotone graph property ${\cal A}$ (a family
of graphs closed under isomorphism and the addition of edges) as:
$$\tau({\cal A}) = \min ~\{0 \leq t \leq \binom{n}{2} :
\widetilde{G}(t) \in {\cal A} \} ~.$$ We use the abbreviation
$\tau(\delta = d)$ for the hitting time $\tau(\{G : \delta(G) \geq
d\})$ of a given graph process, where $\delta(G)$ denotes the
minimal degree of $G$. Finally, we say that a random graph $G$
satisfies some property \textit{with high probability}, or
\textit{almost surely}, or that \textit{almost every} graph process
satisfies a property, if the probability for the corresponding event
tends to $1$ as the number of vertices tends to infinity.

Consider the beginning of the random graph process. It is easy to
see that for every graph $G$, $i(G)$ is at most $\delta(G)$, the
minimal degree of $G$ (choose a set $S$ consisting of a single
vertex of degree $\delta(G)$). Hence, at the beginning of the graph
process, $i(\widetilde{G}(0))=0=\delta(\widetilde{G}(0))$, and this
remains the case as long as there exists an isolated vertex in
$\widetilde{G}(t)$. Next, consider the time where the minimal degree
and maximal degree of the random graph process become more or less
equal. At this point, we can examine random $\delta$-regular graphs
for intuition as to the behavior of the isoperimetric constant, in
which case the results of \cite{AlonUpperRegulerBound} and
\cite{BollobasLowerRegulerBound} imply that $i(\widetilde{G}(t))$ is
roughly $\delta/2$. Hence, at some point along the random graph
process, the behavior of the isoperimetric constant changes, and
instead of being equal to $\delta$ it drifts towards $\delta/2$ (it
is easy to confirm that the isoperimetric constant of the complete
graph is $\frac{n-1}{2}$). The following results summarize the
behavior of the isoperimetric constant of the random graph process
(and, resulting from which, of the appropriate random graphs
models):

In Section \ref{sec-1} we prove that, for almost every graph
process, there is equality between the isoperimetric constant and
the minimal degree, as long as the minimal degree is $o(\log n)$. In
other words, we prove a hitting time result: the minimal degree
increases by 1 exactly when the isoperimetric constant increases by
1 throughout the entire period in which $\delta = o(\log n)$.
\begin{theorem}\label{thm-1}
Let $\ell = \ell(n)$ denote a function satisfying $\ell(n) = o(\log
n)$. Almost every graph process $\widetilde{G}$ on $n$ vertices
satisfies $ i(\widetilde{G}(t)) = \delta(\widetilde{G}(t))$ for
every $t \in [0,\tau(\delta=\ell)]$. Furthermore, with high
probability, for every such $t$, every set $S$ which attains the
minimum of $i(\widetilde{G}(t))$ is an independent set of vertices
of degree $\delta(\widetilde{G}(t))$.
\end{theorem}

In Section \ref{sec-2} we show that the $o(\log n)$ bound in Theorem
\ref{thm-1} is essentially best possible. Indeed, during the period
in which the minimal degree is $\Theta(\log n)$, $i(G)$ drifts
towards $\frac{1}{2}\delta(G)$, as the next theorem demonstrates:
\begin{theorem}\label{thm-2}
For every $0 < \epsilon < \frac{1}{2}$ there exists a constant $C =
C(\epsilon) > 0$, such that the random graph $G \sim
\mathcal{G}(n,p)$, where $p = C \frac{\log n}{n}$, almost surely
satisfies:
$$i(G) \leq \left(\frac{1}{2} + \epsilon\right) \delta(G) = \Theta(\log n)~.$$
Furthermore, with high probability, every set $S$ of size $\lfloor
\frac{n}{2} \rfloor$ satisfies: $\frac{|\partial S|}{|S|} <
\left(\frac{1}{2}+\epsilon\right)\delta(G)$.
\end{theorem}
An analogous statement holds for $\mathcal{G}(n,M)$ as well, where
$M=C n \log n$ for a sufficiently large $C = C(\epsilon)$.

We note that throughout the paper, all logarithms are in the natural
basis.

\section{The behavior of $i(G)$ when $\delta=o(\log n)$}\label{sec-1}
\subsection{Proof of Theorem \ref{thm-1}}
Since every graph $G$ satisfies $i(G) \leq \delta(G)$, proving that,
for every $d \leq \ell$, with high probability, at time $\tau(\delta
= d)$ the isoperimetric constant of $G$ is at least $d$, will prove
the theorem. We show that for every $d = d(n) = o(\log n)$, the
probability for this event is at least $1-o(\frac{1 }{\log n})$, and
the theorem follows from the union bound on the complement events.

Recall that almost every graph process $\widetilde{G}$ satisfies
$\delta(\widetilde{G}) \leq d - 1$ at time $$ m_d = \binom{n}{2}
\frac{ \log n + (d-1)\log\log n- \omega(n)}{n} ~,$$ and
$\delta(\widetilde{G}) \geq d$ at time
$$
M_d = \binom{n}{2} \frac{ \log n + (d-1)\log\log n + \omega(n)}{n}
~,
$$ where $d \geq 1$ is some fixed integer, the $\omega(n)$-term
represents a function growing to infinity arbitrarily slowly while
satisfying $\omega(n) \leq \log\log\log n$, and all logarithms are
natural (see, e.g., \cite{RandomGraphs}, Chapter 3). Hence, $\tau(
\delta=d)$ is between $m_d$ and $M_d$. Using the same methods
described in \cite{RandomGraphs}, it is easy to extend this
statement typically to every $d = d(n) = o(\log n)$, as the next
proposition summarizes:
\begin{proposition}\label{thresholds-prop} Let $\ell = \ell(n) = o(\log n)$. For every $1
\leq d \leq \ell$ define:
$$ r = r(n) = \frac{\log n}{d}~. $$
Next, define the following threshold functions:
\begin{equation}\label{lower-delta-r-bound}
m_d = \binom{n}{2} \frac{ \log n + (d-1)\log r- (2d+\omega(n))}{n}~,
\end{equation}
and:
\begin{equation}\label{upper-delta-r-bound}
M_d = \binom{n}{2} \frac{ \log n + (d-1)\log r + (2d +
\omega(n))}{n} ~,
\end{equation}
where $\omega(n) \leq \log\log r$ and $\lim_{n\rightarrow
\infty}\omega(n)=\infty$. Then, almost every graph process
$\widetilde{G}$ satisfies $\delta(\widetilde{G}(m_d)) \leq d-1$ and
$ \delta(\widetilde{G}(M_d)) \geq d$ for every $1 \leq d \leq \ell$.
\end{proposition}

Notice that $r \leq \log n$, and that $r$ tends to infinity as $n
\rightarrow \infty$, hence these definitions coincide with the
previous definitions of $m_d$ and $M_d$ for a fixed $d$, and it is
left to verify them for $1 \ll d \ll \log n$. Proposition
\ref{thresholds-prop} follows from standard first moment and second
moment considerations, and we postpone its proof to Section
\ref{sec-threshold-prop}. Assume therefore, throughout the proof of
Theorem \ref{thm-1}, that the hitting time $\tau(\delta=d)$ is
almost surely in the interval $(m_d,M_d]$ for every $1 \leq d \leq
\ell$.

Consider a set $S \subset V$ of size $|S|\leq n/2$; we need to show
that, with high probability, every such set satisfies $|\partial S|
\geq \delta(\widetilde{G}(t))|S|$ at every time $t \leq
\tau(\delta=\ell)$ in the random graph process. Clearly, at a given
time $t=M$, the random variable $|\partial S|$ has a binomial
distribution with parameters ${\cal B}\left(|S|(n-|S|),p\right)$,
where $p=M/\binom{n}{2}$. When $|S|$ is sufficiently large (namely,
larger than $n^{1/4}$), the result follows from standard large
deviation bounds and bounds on the tail of the binomial
distribution. However, these bounds are not tight enough for small
values of $|S|$, which require a separate and more delicate
treatment.

Throughout the rest of this section, fix $d = d(n) = o(\log n)$, and
define $m_d$, $M_d$ and $r$ according to Proposition
\ref{thresholds-prop}.

The following lemma shows that every small set $S$ has a boundary of
size at least $\delta(G)|S|$ almost surely:
\begin{lemma}\label{small-sets-lemma}
With probability at least $1-o(n^{-1/5})$, the random graph process
$\widetilde{G}$ satisfies that every $G \in \{\widetilde{G}(t) ~:~
m_d \leq t \leq \tau(\delta = d)\}$ has the property $|\partial S|
\geq \delta(G) |S|$ for every set $S$ of size $|S| \leq n^{1/4}$.
Furthermore, if such a set $S$ satisfies $|\partial S| =
\delta(G)|S|$, it is necessarily an independent set of vertices
whose degrees are $\delta(G)$.
\end{lemma}
\begin{proof}
Given a graph $G=(V,E)$, we call a set $S\subset V$ \textit{bad} if
it satisfies $|\partial S|<\delta(G)|S|$. The idea of the proof is
as follows: we show that, with high probability, every induced
subgraph on $k \leq n^{1/4}$ vertices has a low average degree.
Since bad sets have a boundary of at most $\delta(G) |S|$, this
implies that bad sets, as well as sets which are "almost" bad, must
contain many vertices whose degrees are low in $G$. The result is
derived from several properties of the set of all vertices of low
degrees. We begin with defining this set of vertices and examining
its properties:
\begin{definition} Let $G = (V,E)$. The set of vertices $\Small(G)$ is defined to be:
$$ \Small = \Small(G) = \{v \in V~:~d(v)<4(d+6)\} ~.$$
\end{definition}

\begin{claim}\label{Small-properties}
With probability at least $1-o(n^{-1/5})$, the random graph process
$\widetilde{G}$ has the following property: for every $m_d \leq t
\leq M_d$, $\Small$ is an independent set, and every two vertices of
$\Small$ have no common neighbors in $V$.
\end{claim}
\begin{proof}
Notice that the set $\Small$ changes along the random graph process,
as vertices are removed from it once they reach a degree of
$4(d+6)$. We show a slightly stronger result: if $S_0$ denotes
$\Small(\widetilde{G}(m_d))$, then $S_0$ satisfies the above
properties almost surely for every $m_d \leq t \leq M_d$. Since
$\Small(\widetilde{G}(t)) \subseteq S_0$ for every $t \geq m_d$,
this will imply the claim. In order to prove this result, we show
that, with high probability, $S_0$ satisfies the above properties at
time $t= m_d$, and that the addition of $M_d - m_d$ edges almost
surely does not harm these properties of $S_0$.

Let $p = m_d/\binom{n}{2}$, and let $G_0 \sim \mathcal{G}(n,p)$. The
same consideration will show that $\Small$ satisfies the properties
of the claim with the mentioned probability, both in
$\mathcal{G}(n,p)$ and in $\mathcal{G}(n,m_d)$; for the sake of
simplicity, we perform the calculations in the $\mathcal{G}(n,p)$
model, and note that they hold for the $\mathcal{G}(n,m_d)$ model as
well. Indeed, the main tool in the proof is an upper bound on the
probability for a low degree (a degree of $L=o(n)$ when the edge
probability is $p$), and the probabilities of the relevant events in
$\mathcal{G}(n,m_d)$ are already upper bounded by the corresponding
probabilities in $\mathcal{G}(n,p)$.

Both of the properties mentioned in the claim are immediate
consequences of the next upper bound for the probability of the
event $\{ {\cal B}(n-L,p) \leq D\}$, where $4d \leq D \leq 30d$ and
$L=o(n)$. We use the fact that, by this choice of parameters, $D \ll
(n-L)p$, implying the following monotonicity of the binomial
distribution:
$$\Pr[{\cal B}(n-L,p) \leq D] \leq (D+1)
\binom{n-L}{D}p^{D}(1-p)^{n-L-D} \leq $$
$$ \leq (D+1) \left(\frac{\mathrm{e} p n}{D}\right)^{D} \mathrm{e}^{-p(1-o(1))n} \leq
(30d+1) \left(\frac{(\mathrm{e}+o(1))\log n}{4d}\right)^{30d}
\mathrm{e}^{-(1-o(1))\log n} \leq $$
$$ \leq (30d+1) r^{30d} \mathrm{e}^{-(1-o(1))\log n}
 = \exp\left( O(1) + \log d + 30d \log r -(1-o(1))\log n\right) = $$
 $$ = \exp\left(O(1) + \log d + 30 \log n \frac{\log r}{r} - (1-o(1))\log n\right)
  = \exp\left(-(1-o(1))\log n\right) = o(n^{-0.9})~.$$ Set $D =
4(d+6)$, and let $A_{u,v}$ denote the event that the edge $(u,v)$
belongs to the induced graph on $\Small$, for a given pair of
vertices $u,v \in V$. The following holds:
$$\Pr[A_{u,v}]= p \Pr[{\cal B}(n-2,p) < D-1]^2
\leq \frac{(1+o(1))\log n}{n^{2.8}} = o(n^{-2.5})~.$$ Thus, the
probability that there exists such a pair of vertices is at most
$\binom{n}{2}\Pr[A_{u,v}] = o(n^{-1/2})$, and $\Small(G_0)$ is an
independent set with probability $1-o(n^{-1/2})$. Next, let
$A_{u,v,w}$ denote the event that $u,v \in \Small(G_0)$ and $w$ is a
common neighbor of $u$ and $v$, for some $u,v,w \in V$. Again, we
get:
$$\Pr[A_{u,v,w}] = p^2 \left(p
\Pr[{\cal B}(n-3,p) < D-2]^2 + (1-p) \Pr[{\cal B}(n-3,p) <
D-1]^2\right) \leq $$ $$ \leq p^2 n^{-1.8} = o(n^{-3.5})~,$$ and
therefore $\binom{n}{3}\Pr[A_{u,v,w}] = o(n^{-1/2})$.

We have shown that with probability at least $1-o(n^{-1/2})$,
$\Small(G_0)$ satisfies the two properties of the claim, and by the
same argument, $S_0 = \Small(\widetilde{G}(m_d))$ satisfies the two
properties of the claim with probability at least $1-o(n^{-1/2})$.
We now give a rough upper bound on the size of $S_0$ using the above
upper bound on ${\cal B}(n,p)$:
$$\mathbb{E}|S_0| \leq n \Pr[{\cal B}(n-1,p) < D] = o(n^{0.1})~.$$
Hence, by Markov's inequality, $\Pr[|S_0| \geq n^{0.3}] \leq
n^{-1/5}$. Altogether, we have shown that, with probability
$1-o(n^{-1/5})$, the set $\Small$ at time $t=m_d$ satisfies the
requirements of the claim, and is of size at most $n^{0.3}$.

Assume that indeed $|S_0| \leq n^{0.3}$ and that the distance
between every pair of vertices of $S_0$ is at least $3$ at time
$m_d$. We wish to show that this property is maintained throughout
the period $t \in (m_d,M_d]$. Notice that the probability that an
edge will be added between a given pair of vertices $u,v$ in this
period is
$$\hat{p} = (1+o(1))\left(M_d-m_d\right)/\binom{n}{2} =
(2+o(1))\frac{2d+\omega(n)}{n}~.$$ Hence, the probability that an
internal edge is added to $S_0$ is at most:
$$ \binom{|S_0|}{2} \hat{p} \leq \frac{n^{0.6}(1+o(1))(2d+\omega(n))}{n}
 = o(n^{-1/5}) ~.$$
Since the set of neighbors of $S_0$, $N(S_0)$, consists of at most
$4(d+6)|S_0|$ vertices, the probability that an edge is added
between $N(S_0)$ and a vertex of $S_0$ is at most:
$$ |N(S_0)||S_0|\hat{p} \leq \frac{n^{0.6} (2+o(1)) 4(d+6)(2d+\omega(n))}{n}
= o(n^{-1/5}) ~.$$ Finally, the probability that two edges are added
between one vertex of $V \setminus S_0$ and two vertices of $S_0$ is
at most:
$$ n \binom{|S_0|}{2} \hat{p}^2 \leq \frac{n^{1.6}(2+o(1))
(2d+\omega(n))^2}{n^2}
 = o(n^{-1/5}) ~.$$

Altogether, with probability $1-o(n^{-1/5})$ the set $S_0$ maintains
the property that the distance between each pair of its vertices is
at least $3$ in the period $m_d \leq t \leq M_d$. This completes the
proof of the claim.
\end{proof}

The following claim is crucial to the handling of small sets in $G$,
showing that the average degree in the induced subgraphs on them is
small:
\begin{claim}\label{small-induced-average-degree}
With probability at least $1-o(n^{-1/5})$, the random graph process
$\widetilde{G}$ has the following property: for every $t \leq M_d$,
every induced subgraph of $\widetilde{G}(t)$ on $k \leq n^{1/4}$
vertices contains at most $2k$ edges.
\end{claim}
\begin{proof}
Since this property is monotone with respect to the removal of
edges, it is enough to prove the claim for $t=M_d$. Let
$p=M_d/\binom{n}{2}$ and $G \sim \mathcal{G}(n,p)$. Fix $1 \leq k
\leq n^{1/4}$; the probability that an induced subgraph $H$ on $k$
vertices has at least $2k$ edges is:
$$\Pr[|E(H)| \geq 2k] = \Pr[{\cal B}(\binom{k}{2},p) \geq 2k] \leq
\binom{\binom{k}{2}}{2k} p^{2k} \leq (k p)^{2k} \leq
\left(\frac{(1+o(1))\log n}{n^{3/4}}\right)^{2k} ~.
$$
Summing over all the subgraphs of size at most $k$, we obtain that
the probability that such a subgraph exists is at most: $$\sum_{k
\leq n^{1/4}} \sum_{|H|= k}\Pr[|E(H)| \geq 2k] \leq \sum_{k\leq
n^{1/4}} \binom{n}{k}\left(\frac{(1+o(1)\log n}{n^{3/4}}\right)^{2k}
\leq \sum_{k\leq n^{1/4}} \left(n^{-\frac{1}{2}+o(1)}\right)^{k} =
o(n^{-1/5}) ~.$$ Again, performing the same calculation in
$\mathcal{G}(n,M_d)$ gives the same result: the probability that a
specific set of $2k$ edges belongs to $\mathcal{G}(n,M_d)$ is
$\binom{N-2k}{M_d-2k}/\binom{N}{M_d}$ (where $N=\binom{n}{2}$),
which equals
$\left((1+o(1))M_d/N\right)^{2k}=\left((1+o(1))p\right)^{2k}$.
\end{proof}

Equipped with Claim \ref{Small-properties} and Claim
\ref{small-induced-average-degree}, we are ready to prove Lemma
\ref{small-sets-lemma}.

Recall that a set $S$ is bad iff  $|\partial S|<\delta(G)|S|$. We
call a bad set $S$ \textit{elementary} if it does not contain a
smaller bad set, i.e., every $T \subset S$, $T \neq S$ is not bad.
Clearly, in order to show that there are no bad sets of size at most
$n^{1/4}$, it is enough to show that there are no elementary bad
sets of such size. With high probability, every $G \in
\{\widetilde{G}(t) ~:~ m_d \leq t \leq M_d\}$ satisfies both Claim
\ref{Small-properties} and Claim \ref{small-induced-average-degree}.
Since $m_d < \tau(\delta = d) \leq M_d$, every graph
$G=\widetilde{G}(t)$ in the interval $m_d \leq t \leq
\tau(\delta=d)$ satisfies both claims, as well as $\delta(G) \leq
d$. We claim that this implies the required result; to see this,
consider a graph $G$ which satisfies the above properties, and let
$\delta = \delta(G)$. We first prove that there are no elementary
bad sets of size at most $n^{1/4}$ in $G$:

Let $S$ denote an elementary bad set $S$ of size $k \leq n^{1/4}$.
Notice that necessarily $k\geq 2$, since a single vertex has at
least $\delta$ edges and hence cannot account for a bad set. By
Claim \ref{small-induced-average-degree}, the induced graph $H$ on
$S$ contains at most $2k$ edges. Since the boundary of $S$ contains
at most $\delta k-1 \leq d k$ edges, this implies that $|S \cap
\Small| \geq \frac{3}{4}k$ , otherwise the number of edges in $H$
would satisfy:
$$|E(H)| = \frac{1}{2} \sum_{v \in S}d_H(v) \geq
\frac{1}{2}\left(\frac{k}{4} 4(d+6) - d k \right) \geq 3k~,
$$
leading to a contradiction. Assume therefore that at most $k/4$
vertices in $S$ do not belong to $\Small$. We define $A = S \cap
\Small$, and $B = S \setminus A$. By Claim \ref{Small-properties},
$A$ is an independent set, and furthermore, no two vertices of $A$
have a common neighbor in $B$. Hence, each vertex of $B$ is adjacent
to at most one vertex of $A$, and if we denote by $A' \subseteq A$
the vertices of $A$, which are not adjacent to any vertex of $S$,
the following holds:
$$|A'| \geq |A|-|B| \geq (\frac{3}{4}-\frac{1}{4})k = \frac{1}{2}k~.$$
In particular, $A'$ is nonempty; we claim that this contradicts the
fact that $S$ is elementary. Indeed, each vertex $v \in A'$ is not
adjacent to any vertex in $S$, hence it contributes $d(v)$ edges to
$\partial S$. Removing the vertex $v$ would result in a nonempty ($k
\geq 2$) strictly smaller subset $T$ of $S$ which satisfies:
$$|\partial T| = |\partial S|
- d(v) \leq |\partial S| - \delta < \delta(|S|-1) = \delta|T|~,$$
establishing a contradiction. We conclude that $G$ does not contain
bad sets of size at most $n^{1/4}$.

Next, consider a set $S$ of size $|S| \leq n^{1/4}$ which satisfies
$|\partial S|=\delta |S|$. If $|S|=1$, obviously $S$ consists of a
single vertex of degree $\delta$ and we are done. Otherwise,
repeating the above arguments for bad sets, we deduce that $|S \cap
\Small| \geq \frac{3}{4}|S|$ (this argument merely required that
$|\partial S| \leq \delta|S|$) and that $S$ contains a nonempty set
$A'$, whose vertices are not adjacent to any vertex of $S$. Consider
a vertex $v \in A'$; this vertex contributes $d(v) \geq \delta$
edges to $\partial S$. However, $d(v)$ cannot be greater than
$\delta$, otherwise the set $S' = S \setminus \{v\}$ would satisfy
$|\partial S'| < \delta|S'|$, contradicting the fact that there are
no bad sets of size at most $n^{1/4}$ in $G$. Therefore, all the
vertices of $A'$ are of degree $\delta$, and are not adjacent to any
of the vertices of $S$. If we denote the remaining vertices by $S' =
S \setminus A'$, $S'$ satisfies $|\partial S'| = \delta |S| - \delta
|A'| = \delta |S'|$, and, by induction, the result follows.

This completes the proof of Lemma \ref{small-sets-lemma}.
\end{proof}

The large sets are handled by the following lemma, which shows that
even at time $m_d$ (when the minimal degree is still at most $d-1$)
these sets already have boundaries of size at least $d|S|+1$.
\begin{lemma}\label{large-sets-lemma}
With probability at least $1-o(n^{-1/5})$, the graph
$\widetilde{G}(m_d)$ satisfies $|\partial S|
> d |S|$ for every set $S$ of size $n^{1/4} \leq |S| \leq n/2$ (and
hence $\widetilde{G}(t)$ has this property for every $t \geq m_d$
with probability at least $1-o(n^{-1/5})$).
\end{lemma}
\begin{proof}
Define $p = m_d/\binom{n}{2}$. For the sake of simplicity, the
calculations are performed in the $\mathcal{G}(n,p)$ model and we
note that by the same considerations the results apply for the
corresponding $\mathcal{G}(n,m_d)$ model as well. To show that, with
probability $1-o(n^{-1/5})$, the random graph $G \sim {\cal G}(n,p)$
satisfies $|\partial S| > d |S|$ for sets $S$ of the given size,
argue as follows:

Fix a set $S \subset V$ of size $k$, $\frac{n}{\log n} \leq k \leq
n/2$, and let $A_S$ denote the event $\{|\partial S| \leq d k\}$.
Let $\mu$ denote $\mathbb{E}|\partial S| = k (n-k) p$. By the
Chernoff bound, $\Pr[|\partial S| < \mu - t] \leq
\exp\left(-\frac{1}{2\mu}t^2\right)$. Therefore, setting $t=\mu -
(dk+1)$, we get:
$$ \Pr[A_S] = \Pr[|\partial S| < d k + 1] \leq
\exp\left(-\frac{1}{2}\left(1-\frac{d+\frac{1}{k}}{(n-k)p}\right)^2
k(n-k)p\right) \leq $$
$$ \leq \exp\left(-\frac{1}{2}\left(1-\frac{(2+o(1))d}{\log n}\right)^2 k
\left(\frac{1}{2}-o(1)\right)\log n \right) =
\exp\left(-\frac{1-o(1)}{4} k \log n \right) ~.$$ Hence, the
probability that there exists such a set $S$ is at most:
$$ \sum_{k=\frac{n}{\log n}}^{n/2}\binom{n}{k}\exp
\left(-\frac{1-o(1)}{4} k \log n \right) \leq \sum_{k=\frac{n}{\log
n}}^{n/2} \left(\mathrm{e} \frac{n}{k}\right)^k
\exp\left(-\frac{1-o(1)}{4} k \log n \right) \leq $$
$$ \leq \sum_{k=\frac{n}{\log n}}^{n/2} \exp\left( k(\log\log n + 1)
-\frac{1-o(1)}{4} k \log n \right) \leq \sum_{k=\frac{n}{\log
n}}^{n/2} \left(n^{-\frac{1}{4}+o(1)}\right)^k = o(n^{-1/5})~.
$$
Let $S \subset V$ be a set of size $n^{1/4} \leq k \leq
\frac{n}{\log n}$. Notice that:
\begin{equation}\label{d-(n-k)p-eq-(1+o(1))}(n-k)p= (1+o(1))\log n
~,
\end{equation}
and hence, $d k < \mu$, and we can give the following upper bound on
the probability that $|\partial S| \leq d k$:
$$ \Pr[ |\partial S| \leq d k] \leq (d k+1) \Pr[|\partial S| =
d k] = (d k+1) \binom{k(n-k)}{d k} p^{d k}(1-p)^{k(n-k) - d k} \leq
$$
$$ \leq (d k +1) \left(\frac{\mathrm{e} k(n-k)p}{d k}\right)^{d k }
 \mathrm{e}^{-p k (n-k-d)} = (d k +1)
(\mathrm{e}/d)^{d k} \left(p(n-k)\right)^{d k} \mathrm{e}^{-k p n +
p k^2 + p k d}~.$$ We now use \eqref{d-(n-k)p-eq-(1+o(1))} and the
facts that $p k \leq 1+o(1)$ and  $d = o(k)$, and obtain:
$$ \Pr[ |\partial S| \leq d k] \leq
O(1)d k (\mathrm{e}/d)^{d k} \frac{(\log n)^{d k}
\mathrm{e}^{\left(2d+\omega(n)+1+o(1)\right)k+d} }{n^k r^{k(d - 1)}}
\leq \left(\frac{\mathrm{e}^{\omega(n)+2d+O(1)} \log n}{n}\right)^k
~.$$ Summing over all sets $S$ of size $k$, we get:
$$ \sum_{|S|=k} \Pr[|\partial S|\leq d k] \leq
\left(\frac{\mathrm{e}n}{k} \right)^k
\left(\frac{\mathrm{e}^{\omega(n)+2d+O(1)} \log n}{n}\right)^k =
\left(\frac{\mathrm{e}^{\omega(n)+2d+O(1)} \log n}{k}\right)^k
 \leq $$
$$ \leq \left(\frac{O(1)n^{2/r} \log r \log n}{n^{1/4}}\right)^k =
\left(n^{-\frac{1}{4}+o(1)}\right)^k~.$$ Thus:
$$ \sum_{n^{1/4} \leq |S|\leq \frac{n}{\log n}} \Pr[|\partial S| \leq d|S|]
\leq \sum_{k \geq n^{1/4}}  \left(n^{-\frac{1}{4}+o(1)}\right)^k =
o(n^{-1/5})~.
$$
\end{proof}

Since $\widetilde{G}$ satisfies the properties of both Lemma
\ref{small-sets-lemma} and Lemma \ref{large-sets-lemma} for a given
$d \leq \ell = o(\log n)$ with probability at least $1-o(n^{-1/5})$,
the union bound over all possible values of $d$ implies that these
properties are satisfied almost surely for every $d \leq \ell$.
Theorem \ref{thm-1} follows directly: to see this, assume that
indeed a random graph process $\widetilde{G}$ satisfies the
mentioned properties for every $d \leq \ell$, and consider some $d
\leq \ell$. By the properties of Lemma \ref{small-sets-lemma}, in
the period $t\in [m_d,\tau(\delta = d)]$ every set of size $k \leq
n^{1/4}$ has at least $\delta k$ edges in its corresponding cut, and
if there are precisely $\delta k$ edges in the cut, then $S$ is an
independent set of vertices of degree $\delta$. In particular, at
time $t = \tau(\delta = d)$, every set $S$ of at most $n^{1/4}$
vertices has a ratio $\frac{|\partial S|}{|S|}$ of at least $d$, and
a ratio of precisely $d$ implies that $S$ is an independent set of
vertices of degree $d$. By monotonicity, this is true for every $t
\in \left[\tau(\delta=d),\tau(\delta=d+1)\right)$. Next, by the
properties of Lemma \ref{large-sets-lemma}, every set of size $k
\geq n^{1/4}$ has at least $d k + 1$ edges in its corresponding cut
at time $t=m_d$. In particular, for every $t \in
[\tau(\delta=d),\tau(\delta=d+1))$, every set $S$, larger than
$n^{1/4}$ vertices, has a ratio $\frac{|\partial S|}{|S|}$ strictly
larger than $d$. These two facts imply the theorem.
 \qed

\subsection{Proof of Proposition \ref{thresholds-prop}}
\label{sec-threshold-prop} A standard first moment consideration
shows that indeed, with high probability,
$\delta(\mathcal{G}(n,M_d)) \geq d$ for every $d \leq \ell$. We
perform the calculations in the $\mathcal{G}(n,p)$ model and note
that the same applies to $\mathcal{G}(n,M_d)$.

For each $v \in V(G)$, let $A_v$ and $B_v$ denote the events $\{
d(v) = d - 1\}$ and $\{ d(v) \leq d-1\}$ respectively, and set
$Y_d=|\{v:d(v)=d-1\}|$ and $Z_d=|\{v:d(v)\leq d-1\}|$. Recall that
$d = o(\log n)$, and furthermore, we may assume that $d$ tends to
infinity as $n \rightarrow \infty$, since $m_d$ and $M_d$ coincide
with the well known threshold functions for constant values of $d$.
Choosing $p=M_d/\binom{n}{2}$, the following holds:
$$ \Pr[A_v] = \binom{n-1}{d-1}p^{d-1}(1-p)^{n-d} \leq
\left(\frac{(1+o(1))\mathrm{e}\log n}{d}\right)^{d-1}
\mathrm{e}^{-(1-\frac{d}{n})\left(\log n + (d-1)\log r + 2d +
\omega(n)\right)} \leq $$
$$ \leq \frac{1}{n^{1-d/n}} \left( \frac{(1+o(1))\mathrm{e} r}
{r^{1-d/n}}\right)^{d-1}\mathrm{e}^{-(1-o(1))(2d+\omega(n))}=$$
\begin{equation}\label{threshold-A-v-bound-M-d}
= \frac{n^{d/n}}{n} \left((1+o(1))\mathrm{e} r^{\frac{\log n}{r
n}}\right)^{d-1}\mathrm{e}^{-(1-o(1))(2d+\omega(n))} \leq
\frac{1}{n}\mathrm{e}^{-(1-o(1))(d+\omega(n))}~.
\end{equation}
Since $d \leq (n-1)p$, we have:
$$\Pr[B_v] \leq d \Pr[A_v] \leq \frac{1}{n}\mathrm{e}^{-(1-o(1))(d+\omega(n))}~.
$$
Hence,
$$ \mathbb{E}Z_d \leq \mathrm{e}^{-(1-o(1))
(d+\omega(n))}~,
$$
and summing over every $d \leq \ell$ we obtain:
$$\sum_{1 \ll d \leq \ell} \Pr[Z_d>0] \leq
\mathrm{e}^{-(1-o(1))\omega(n)} \sum_{1 \ll d \leq \ell}
\mathrm{e}^{-(1-o(1))d} = o(1) ~.$$

A second moment argument proves that almost surely
$\delta(\mathcal{G}(n,p)) \leq d-1$ for every $d \leq \ell$. To see
this, argue as follows (again, calculations are performed in the
$\mathcal{G}(n,p)$ model): following the same definitions, only this
time with $p=m_d/\binom{n}{2}$, apply the bound $\binom{a}{b} \geq
\left(\frac{a}{b}\right)^b$ and the well known bound $1-x\geq
\mathrm{e}^{-x/(1-x)}$ for $0 \leq x<1$, to obtain:
$$ \Pr[A_v] = \binom{n-1}{d-1}p^{d-1}(1-p)^{n-d} \geq $$
$$ \geq
\left(\frac{(1+o(1))\log n}{d}\right)^{d-1} \mathrm{e}^{\left(-\log
n - (d-1)\log r + 2d +\omega(n)\right)/(1-p)} \geq \frac{1}{n}
\Omega(\mathrm{e}^{d+\omega(n)}) ~,$$ where in the last inequality
we omitted the the $1/(1-p)$ factor in the exponent, since, for
instance, $n^{1-\frac{1}{1-p}}=n^{\frac{-p}{1-p}} \geq
n^{-O(1)\frac{\log n}{n}} = \mathrm{e}^{o(1)}$. Therefore:
$$\mathbb{E}Y_d = \Omega(\mathrm{e}^{d+\omega(n)})~.$$
Take $u \neq v \in V(G)$; denoting $\mathcal{P}^K_L =
\Pr[\mathcal{B}(K,p)=L]$, the following holds:
$$ \textup{Cov}(A_u,A_v) = \Pr[A_u \wedge A_v] - \Pr[A_u] \Pr[A_v] =
p (\mathcal{P}^{n-2}_{d-2})^2 + (1-p)(\mathcal{P}^{n-2}_{d-1})^2 -
(\mathcal{P}^{n-1}_{d-1})^2 ~.$$ Since $\mathcal{P}^{n-1}_{d-1}=p
\mathcal{P}^{n-2}_{d-2} + (1-p)\mathcal{P}^{n-2}_{d-1}$, we get:
\begin{eqnarray} \textup{Cov}(A_u,A_v) &=& p (1-p) (\mathcal{P}^{n-2}_{d-2})^2 +
(1-p)p (\mathcal{P}^{n-2}_{d-1})^2 -2p(1-p)\mathcal{P}^{n-2}_{d-2}
\mathcal{P}^{n-2}_{d-2} = \nonumber \\
 &=& p(1-p)(\mathcal{P}^{n-2}_{d-1}-\mathcal{P}^{n-2}_{d-2})^2 \leq
 p(\mathcal{P}^{n-2}_{d-1})^2 \nonumber ~.\end{eqnarray}
 Notice that $\mathcal{P}^{n-2}_{d-1}$ corresponds to the event $A_v$
 for a graph on $n-1$ vertices, and hence a similar calculation to
 the one in \eqref{threshold-A-v-bound-M-d}
 shows that $\mathcal{P}^{n-2}_{d-1} =
 O(\exp(3d+\omega(n))/n)$.
 Altogether we get:
 $$ \textup{Cov}(A_u,A_v) \leq
 O(1)p \frac{\mathrm{e}^{6d+2\omega(n)}}{n^2} \leq
 O(1) \mathbb{E}Y_d \frac{\mathrm{e}^{5d+\omega(n)}\log n }{n^3}
 = o(n^{-2}) \mathbb{E}Y_d~,$$
which gives the following upper bound on the variance of $Y_d$:
$$ \textup{Var}(Y_d) \leq \mathbb{E}Y_d + \sum_{u \neq v}\textup{Cov}(A_u,A_v)
\leq \mathbb{E}Y_d + n^2 o(n^{-2}) \mathbb{E}Y_d =_d
(1+o(1))\mathbb{E}Y_d~.
$$ Applying Chebyshev's inequality gives:
$$\Pr[Y_d = 0] \leq \frac{\textup{Var}(Y_d)}{(\mathbb{E}Y_d)^2} \leq \frac{1+o(1)}{\mathbb{E}Y_d} \leq
O(\mathrm{e}^{-d-\omega(n)})~, $$ and summing over every $d \leq
\ell$ we obtain:
$$\sum_{1 \ll d \leq \ell} \Pr[Y_d=0] \leq
O(1)\mathrm{e}^{-\omega(n)} \sum_{1 \ll d \leq \ell} \mathrm{e}^{-d}
= o(1) ~,$$ as required. \qed

\section{The behavior of $i(G)$ when $\delta=\Omega(\log n)$}\label{sec-2}
\textbf{Proof of Theorem \ref{thm-2}}: { A bisection of a graph $G$
on $n$ vertices is a partition of the vertices into two disjoint
sets $(S, T)$, where $|S| = \lfloor \frac{n}{2} \rfloor$ and $T = V
\setminus S$. Fix $\epsilon_1 > 0$; we first prove that, with high
probability, every bisection $(S,V\setminus S)$ of
$\mathcal{G}(n,p)$ has strictly less than
$\left(\frac{1}{2}+\epsilon_1\right)n p |S|$ edges in the cut it
defines, provided that $\lim_{n\rightarrow \infty} n p = \infty$. We
omit the floor and ceiling signs to simplify the presentation of the
proof.

Let $S$ be an arbitrary set of $n/2$ vertices. The number of edges
in the boundary of $S$ has a binomial distribution with parameters
${\cal B}(n^2/4,p)$, hence (by our assumption on $p$) its expected
value $\mu$ tends to infinity faster than $n$. By the Chernoff
bound, $\Pr[|\partial S| \geq (1+t)\mu]\leq \exp(-\mu t^2 / 4)$
provided that $t<2\mathrm{e}-1$, thus we get:
$$ \Pr[|\partial S| \geq \left(\frac{1}{2}+\epsilon_1\right) n p |S|]
= \Pr[|\partial S| \geq \left(1 + 2\epsilon_1\right) \mu]
 \leq \exp(-\Omega(\mu))  ~.$$ Since this
probability is $o(2^{-n})$, the expected number of bisections, in
which the corresponding cuts contain at least
$\left(\frac{1}{2}+\epsilon_1\right) n p |S|$ edges, is $o(1)$.

Next, fix $\epsilon_2 > 0$. We claim that the minimal degree of
$\mathcal{G}(n,p)$, where $p=C \frac{\log n}{n}$ and
$C=C(\epsilon_2)$ is sufficiently large, is at least
$(1-\epsilon_2)n p$. Applying the Chernoff bound on the binomial
distribution representing the degree of a given vertex $v$ gives:
$$\Pr[d(v) \leq (1-\epsilon_2) n p] = \Pr[d(v) \leq (1-\epsilon_2+o(1))
\mathbb{E}d(v)] \leq \exp\left(-C \frac{\epsilon_2^2}{2}(1-o(1))
\log n\right)~,$$ and for $C > \frac{2}{\epsilon_2^2}$ this
probability is smaller than $\frac{1}{n}$.

Altogether, for a sufficiently large $C$, the following holds with
high probability: every bisection $(S,V\setminus S)$ satisfies:
$$
\frac{|\partial S|}{|S|} <
\frac{\frac{1}{2}+\epsilon_1}{1-\epsilon_2}\delta(G) =
\left(\frac{1}{2}+\epsilon\right) \delta(G)~,$$ where
$\epsilon=\frac{\epsilon_1+\epsilon_2/2}{1-\epsilon_2}$. \qed}

{\remark We note that the above argument gives a crude estimate on
the value of $C = C(\epsilon)$. Since the first claim, concerning
the behavior of bisections, holds for every value of $C$, we are
left with determining when typically the minimal degree of $G$
becomes sufficiently close to the average degree. This threshold can
be easily computed, following arguments similar to the ones in the
proof of Proposition \ref{thresholds-prop}; the following value of
$C(\epsilon)$ is sufficient for the properties of the theorem to
hold with high probability:
$$
C > \frac{1+2\epsilon}{2\epsilon-\log (1+2\epsilon)} ~.$$ }

{\remark Theorem \ref{thm-2} provides an upper bound on $i(G)$,
which is almost surely arbitrarily close to $\frac{\delta}{2}$ while
the graph satisfies $\delta=\Theta(\log n)$. We note that the
arguments of Theorem \ref{thm-1} can be repeated (in a simpler
manner) to show that with high probability $i(\widetilde{G}(t)) \geq
\delta/2$ for every $t$, and hence the bound in Theorem \ref{thm-2}
is tight.}

\section{Concluding remarks}

We have shown that there is a phase transition when the minimal
degree changes from $o(\log n)$ to $\Omega(\log n)$; it would be
interesting to give a more accurate description of this phase
transition. Theorem \ref{thm-1} treats $\delta(G)=o(\log n)$, and
Theorem \ref{thm-2} shows that, almost surely, $i(G) < \delta(G)$
once $p = C \log n/n$, where $X > 2/(1-\log2) \approx 6.52$, in
which case $\delta(G) > (C/2) \log n$. Hence we are left with the
period in which $\delta(G) = c \log n$, where $0 < c \leq
1/(1-\log2) \approx 3.26$. It seems plausible to show that in this
period $i(G) = \delta(G)$, i.e., that the isoperimetric constant is
determined either by the typical minimal degree, or by the typical
size of a bisection.

The vertex version of the isoperimetric constant (minimizing the
ratio $|\delta S|/|S|$, where $\delta S \subset V \setminus S$ is
the vertex neighborhood of $S$) is less natural, since the minimum
has to be defined on all nonempty sets of size at most
$n/(K+\epsilon)$ if we wish to allow the constant to reach the value
$K$. Nevertheless, the methods used to prove Theorem \ref{thm-1} can
prove similar results for the vertex case, at least as long as the
minimum degree is constant. Indeed, in that case, the probability
for two vertices to have a common neighbor is small enough not to
have an effect on the results.

Finally, it is interesting to consider the isoperimetric constant of
certain subgraphs along the random graph process. To demonstrate
this, we consider the period of $\widetilde{G}$ in which the minimal
degree is $0$, i.e., $t \leq \tau(\delta=1)$. The existence of
isolated vertices in $\widetilde{G}(t)$ implies that
$i(\widetilde{G}(t))=0$, however even if we disregard these
vertices, and examine $G'(t)$, the induced subgraph on the
non-isolated vertices, then after a short while (say, at $t = c n$
for some $c>0$), $i(G'(t)) < \epsilon$ for every $\epsilon
> 0$. An easy calculation shows that small sets, with high
probability, have an edge boundary which is smaller than their size.
For instance, when $p = c/n$ for some $c < 1$, $\mathcal{G}(n,p)$
almost surely satisfies that all connected components are of size
$O(\log n)$, hence each component $\mathcal{C}$ has a ratio
$\frac{|\partial \mathcal{C}|}{|\mathcal{C}|}$ of $0$. Furthermore,
if we take $p = C/n$ for some $C>1$, and consider the giant
component $H$ (recall that for this value of $p$, almost surely
there is a single component of size $\Theta(n)$, and all other
components are of size $O(\log n)$), $i(H)<\epsilon$ for every
$\epsilon>0$. One way to see this, is to consider a collection of
arbitrarily long paths, each of which connects to the giant
component at precisely one end.

\noindent\textbf{Acknowledgement} The authors would like to thank
Noga Alon for helpful discussions and keen observations.


\begin{thebibliography}{alpha}

\bibitem{AlonUpperRegulerBound}
N. Alon, On the edge-expansion of graphs, Combinatorics, Probability
and Computing 6 (1997), 145-152.

\bibitem{AlonLowerEigenvalueBound}
N. Alon and V.D. Milman, $\lambda_1$, isoperimetric inequalities for
graphs and superconcentrators, J. Combinatorial Theory, Ser. B 38
(1985), 73-88.

\bibitem{BollobasLowerRegulerBound}
B. Bollob\'{a}s, The isoperimetric number of random regular graphs,
Europ. J. Combinatorics 9 (1988), 241-244.

\bibitem{BollobasLeaderTorus}
B. Bollob\'{a}s and I. Leader, An isoperimetric inequality on the
discrete torus, SIAM J. Disc. Math. 3 (1990) 32-37.

\bibitem{BollobasLeaderGrid}
B. Bollob\'{a}s and I. Leader, Edge-isoperimetric inequalities in
the grid, Combinatorica 11(1991) 299-314.

\bibitem{RandomGraphs}
B. Bollob\'{a}s, Random graphs, volume 73 of Cambridge Studies in
Advanced Mathematics, Cambridge University Press, Cambridge, second
edition, 2001.

\bibitem{Buser}
P. Buser, Cubic graphs and the first eigenvalue of a Riemann
surface, Mathematische Zeitschrift, 162 (1978), 87-99.

\bibitem{ChungCheeger}
F.R.K. Chung, Laplacians of graphs and Cheeger's inequalities, in:
Proc. Int. Conf. "Combinatorics, Paul Erd\H{o}s is Eighty",
Keszthely (Hungary), 1993, 2, 1–16.

\bibitem{SpectralGraphTheory}
F.R.K. Chung, Spectral Graph Theory, American Mathematical Society,
no. 92 in the Regional Conference Series in Mathematics, Providence,
RI, 1997.

\bibitem{ChungTetali}
F.R.K. Chung and P. Tetali, Isoperimetric inequalities for cartesian
products of graphs, Combinatorics, Probability and Computing 7
(1998), 141-148.

\bibitem{FriedmanKahnSzemeredi}
J. Friedman, J. Kahn and E. Szemer\'{e}di, On the second eigenvalue
in random regular graphs, Proc. 21st ACM STOC (1989), 587-598.

\bibitem{HoudreTetali}
C. Houdr\'{e} and P. Tetali, Isoperimetric constants for product
Markov chains and graph products, Combinatorica Vol. 24 (2004),
359-388.

\bibitem{MoharIso}
B. Mohar, Isoperimetric numbers of graphs, Journal of Combinatorial
Theory, Series B, 47 (1989), 274-291.

\bibitem{ShamirUpfal}
E. Shamir and E. Upfal, On factors in random graphs, Israel J. Math.
39 (1981), no. 4, 296-302.

\end{thebibliography}
\end{document}